%% file: qgor-new.tex
\documentclass{amsart}

\usepackage{amsmath}
\usepackage{amssymb}
\usepackage{amscd}
\usepackage{amsthm}
\usepackage{enumerate}
\usepackage{array}
\usepackage{epsf} 
\usepackage{pictexwd}
\usepackage{rotating}
\usepackage{float}  \restylefloat{figure} 
\usepackage[OT2,OT1]{fontenc}

\input macros

\begin{document}

\input topmatter

\maketitle

\input main

\bibliographystyle{amsalpha}
\bibliography{qgor-new}

\end{document}

%% file: macros.tex
\theoremstyle{plain}
\newtheorem*{thmA}{Theorem A}
\newtheorem*{thmB}{Theorem B}
\newtheorem{lem}{Lemma}
\newtheorem{prop}[lem]{Proposition}

\theoremstyle{definition}

\theoremstyle{remark}
\newtheorem*{rem}{Remark}

\newdimen\XX

\def\su{\mathfrak{su}}

\def\forcehmode{\hskip0pt\relax}

\let\myskip=\medskip

\font\msbm=msbm10
\def\semiprod{\hbox{\msbm\char111}}

\def\st{\,\,\big|\,\,}
\def\bs{\backslash}
\def\<{\langle}
\def\>{\rangle}
\def\ie{i.e.\ }

\let\ge=\geqslant
\let\le=\leqslant
\let\emptyset=\varnothing

\def\definebb#1=#2.{\def#1{{{\mathbb #2}^{\vphantom{x}}}}}
\definebb \c=C.
\definebb \r=R.
\definebb \d=D.
\definebb \n=N.
\definebb \z=Z.
\definebb \q=Q.

\def\calF{\mathcal F}

\DeclareMathOperator{\Cl}{Cl}

\DeclareMathOperator{\Id}{Id}
\DeclareMathOperator{\Int}{Int}
\DeclareMathOperator{\Isom}{Isom}
\DeclareMathOperator{\lcm}{lcm}

\DeclareMathOperator{\MathOpPSU}{PSU}
\DeclareMathOperator{\MathOpSU}{SU}

\let\Re=\undefined \DeclareMathOperator{\Re}{Re}

\def\dd{\partial}
\def\ddP{{\dd P}}

\def\al{{\alpha}}
\def\ve{{\varepsilon}}
\def\thet{{\vartheta}}
\def\De{{\Delta}}

\def\PSU{\MathOpPSU(1,1)}
\def\SU{\MathOpSU(1,1)}
\def\tSU{\widetilde{\MathOpSU}(1,1)}

\def\G{{\Gamma}}
\def\bG{\overline\G}
\def\la{{\lambda}}

\def\Grp{G}
\def\tGrp{\tilde\Grp}
\def\Lrp{L}
\def\tLrp{\tilde\Lrp}

\def\bdlmap{\theta}
\def\covbdlmap{\theta}
\def\covmap{\pi}

\def\WZZW{\begin{pmatrix} w&z \\ \bar z&\bar w \end{pmatrix}}

\def\Eg{E_g}
\def\Ig{I_g}
\def\Hg{H_g}

\def\Fg{F_g}
\def\Eh{E_h}
\def\Fh{F_h}
\def\Ia{I_a}
\def\Ee{E_e}
\def\Ie{I_e}
\def\He{H_e}
\def\Xe{X_e}
\def\Fe{F_e}

\def\calFg{\calF_g}

\def\barg{{\bar g}}
\def\bare{{\bar e}}

\def\bE{{\bar E}}
\def\bI{{\bar I}}

\def\bEbg{\bE_\barg}
\def\bIbg{\bI_\barg}
\def\bEbe{\bE_\bare}
\def\bIbe{\bI_\bare}

\def\Tx{T(x)}
\def\Tu{T(u)}
\def\Qx{Q_x}
\def\Qu{Q_u}
\def\Pu{P_u}
\def\Xu{X_u}

\def\deck{d}

\def\spieg{\eta}
\def\inv{\ve}

\def\Rx{R_x}

\def\Gu{{\G(u)}}
\def\Guou{{\Gu\backslash\{u\}}}

\def\cupddQxGu{{\bigcup\limits_{\hbox to 0pt{\hss$\scriptstyle x\in\Gu$\hss}}\dd\Qx}}
\def\capRxGu{{\bigcap\limits_{\hbox to 0pt{\hss$\scriptstyle x\in\Gu$\hss}}\Rx}}
\def\cupQxGu{{\bigcup\limits_{\hbox to 0pt{\hss$\scriptstyle x\in\Gu$\hss}}\Qx}}
\def\capRxGuou{{\bigcap\limits_{\hbox to 0pt{\hss$\scriptstyle x\in\Guou$\hss}}\Rx}}
\def\cupIntQxGuou{{\bigcup\limits_{\hbox to 0pt{\hss$\scriptstyle x\in\Guou$\hss}}\Int\Qx}}

\def\hF{\hat F}

\newenvironment{pictexure}[2]{\bgroup
  \def\PICdirectory{#1}\beginpicture
  \input #1/#2.pic \ignorespaces}{\endpicture\egroup}

\newcommand{\Pictexure}[2]{\begin{pictexure}{#1}{#2}\relax\end{pictexure}}

\let\cdlabelsize=\small
\def\cdlhss{0.5em\relax}
\newdimen\defaultarrow \defaultarrow 28pt
\newdimen\arrskip \arrskip=0.11\defaultarrow
\newif\ifcdarrowdepth \cdarrowdepthtrue
\newif\ifcdarrowheight \cdarrowheighttrue
\def\spacedfill #1#2{
  \hskip#1\arrskip\hbox to#1\defaultarrow{\cdlabelsize#2}\hskip#1\arrskip}

\def\cdarrow#1#2{
   \ifnum#1=0  \spacedfill{#2}{\rightarrowfill}\else
   \ifnum#1=180\spacedfill{#2}{\leftarrowfill }\else
   $\vcenter{\hbox{\vtop{\hrule height0pt
  \hbox{\begin{turn}{#1}%
   \hbox{\spacedfill{#2}{\rightarrowfill}}%
  \end{turn}}}}}$\fi\fi}

\def\cdarrowxx#1#2#3#4{
  \putcdlabel {\cdlabelsize$\displaystyle\genfrac{}{}{0pt}{}{#1}{#2}$}
   onto box {\cdarrow{#3}{#4}}}

\def\cdarrowv#1#2#3#4{
  \llap{\cdlabelsize$#1$\hskip\cdlhss}%
  \cdarrow{#3}{#4}%
  \rlap{\cdlabelsize\hskip\cdlhss$#2$}}

\newbox\arrowbox
\def\putcdlabel #1 onto box #2{
  \setbox\arrowbox=\hbox{#2}%
  \setbox0=\hbox{\rlap{\copy\arrowbox}\hbox to\wd\arrowbox{\hss #1\hss}}%
  \ifcdarrowheight\ht0=\ht\arrowbox\relax\fi \cdarrowheighttrue
  \ifcdarrowdepth\dp0=\dp\arrowbox\relax\fi \cdarrowdepthtrue
  \box0}

%% file: topmatter.tex
\author[Anna Pratoussevitch]{Anna Pratoussevitch}
\address{Department of Mathematical Sciences\\ University of Liverpool\\ Peach Street \\ Liverpool L69~7ZL}
\email{annap@liv.ac.uk}

\title[The Combinatorial Geometry of $\mathbb Q$-Gorenstein Singularities]%
{The Combinatorial Geometry of $\mathbb Q$-Gorenstein Quasi-Homogeneous Surface Singularities}

\begin{date}  {\today} \end{date}


\begin{abstract} \input abstract \end{abstract}

\subjclass[2000]{Primary 53C50; Secondary 14J17, 32S25, 51M20, 52B10}






\keywords{Lorentz space form, polyhedral fundamental domain, quasihomogeneous singularity.}

%% file: abstract.tex

The main result of this paper is a construction of fundamental domains for certain group actions on Lorentz manifolds of constant curvature.
We consider the simply connected Lie group $\tGrp=\widetilde{\operatorname{SU}}(1,1)$.
The Killing form on the Lie group~$\tGrp$ gives rise to a biinvariant Lorentz metric of constant curvature.
We consider a discrete subgroup $\Gamma_1$ and a cyclic discrete subgroup $\Gamma_2$ in~$\tGrp$ which satisfy certain conditions.
We describe the Lorentz space form $\Gamma_1\backslash\tGrp/\Gamma_2$ by constructing a fundamental domain
for the action of $\Gamma_1\times\Gamma_2$ on~$\tGrp$ by~$(g,h)\cdot x=gxh^{-1}$.
This fundamental domain is a polyhedron in the Lorentz manifold $\tGrp$ with totally geodesic faces.
For a co-compact subgroup the corresponding fundamental domain is compact.


%% file: main.tex
\section{Introduction}

In the context of Riemannian manifolds, there are standard constructions for fundamental domains, for example Dirichlet regions.
However, in the context of semi-Riemannian manifolds, such constractions are rare.
The main result of this paper is a construction of fundamental domains for certain group actions on Lorentz manifolds of constant curvature.

\myskip
We consider the universal cover~$\tGrp$ of the group~$\Grp$ of orien\-ta\-tion-preserving isometries of the hyperbolic plane.
The Killing form on the Lie group~$\tGrp$ gives rise to a biinvariant Lorentz metric of constant curvature.
We consider a discrete subgroup~$\G_1$ and a discrete cyclic subgroup~$\G_2$ in~$\tGrp$ which satisfy the conditions~$(*)$ specified below.
In this paper we describe a construction of fundamental domains for the action of $\G_1\times\G_2$ on $\tSU$ by~$(g,h)\cdot x=gxh^{-1}$.
The resulting fundamental domain is a polyhedron in the Lorentz manifold $\tGrp$ with totally geodesic faces.
For a co-compact subgroup the corresponding fundamental domain is compact.
The precise formulation of these results is contained in Theorems A and B.

\myskip
The study of discrete subgroups of finite level is motivated by some deep connections
between these subgroups and quasi-homogeneous isolated singularities of complex surfaces studied by J.~Milnor, I.~Dolgachev, and W.~Neumann~\cite{Mi, Do83, Neu77, Neu83}.
The class of subgroups for which we construct fundamental domains corresponds to an interesting class of singularities.
There is a 1-1-correspondence between the subgroups from this class and quasi-homogeneous $\q$-Gorenstein surface singularities.
In particular the bi-quotients of the form~$\G_1\backslash\tGrp/\G_2$ are diffeomorphic to the links of quasi-homogeneous $\q$-Gorenstein singularities.
For a more detailed treatment of this connection see~\cite{Pr:qgor} and~\cite{BPR}, \S 1--2.

\myskip
The construction described in~\cite{Pr:diss}, \cite{BPR}, \cite{Pr:statia} can be understood as a special case of the construction described in this paper
when the subgroup~$\G_2$ is trivial.


\myskip
A bi-quotient of the form~$\G_1\backslash\tGrp/\G_2$ is a standard Lorentz space form.
The standard Lorentz space forms were studied by R.S.~Kulkarni and F.~Raymond~\cite{KR}.
Examples of non-standard Lorentz space forms were found by W.~Goldman~\cite{Go}, \'E.~Ghys~\cite{Gh}, and recently by F.~Salein~\cite{Sa}.
The survey~\cite{BZ} of Th.~Barbot and A.~Zeghib and the paper~\cite{F} of Ch.~Frances are good references for the reader interested in group actions on Lorentz manifolds.
The results of this paper suggest that the description of Lorentz space forms by means of fundamental domains could be extended to include non-standard Lorentz space forms.

\myskip

\myskip
Let us specify the conditions that we want to impose on the subgroups~$\G_1$ and~$\G_2$.
We consider the universal cover of the group~$\Grp=\PSU$ of orien\-ta\-tion-preserving isometries of the hyperbolic plane.
Here our model of the hyperbolic plane is the unit disc $\d$ in $\c$.
The kernel of the universal covering map $\tSU\to\PSU$ is the centre $Z$ of the group $\tSU$, an infinite cyclic group.
Therefore, for each natural number~$k$ there is a unique connected $k$-fold covering of $\PSU$.
For $k=2$ this is the group
$$\SU=\left\{\WZZW\st(w,z)\in\c^2,~|w|^2-|z|^2=1\right\}.$$
The {\sl level\/} of a discrete subgroup $\G\subset\tSU$ is the index of $\G\cap Z$ as a subgroup of $Z$.

\myskip
{\bf Condition~(*):}
We consider a discrete subgroup~$\G_1$ and a discrete cyclic subgroup~$\G_2$ in~$\tSU$ of finite level~$k$.
We suppose that the images~$\bG_1$ and~$\bG_2$ of $\G_1$ resp.~$\G_2$ in~$\PSU$ have a joint fixed point in $\d$,
\ie there is a point~$u$ in~$\d$ which is fixed by a nontrivial element of~$\bG_1$ and by a nontrivial element of~$\bG_2$.
For~$i=1,2$, let $p_i$ be the smallest order of a non-trivial element in~$\bG_i$ that has~$u$ as a fixed point.
Let $p=\lcm(p_1,p_2)$ be the least common multiple of~$p_1$ and~$p_2$.
Furthermore we assume that~$p>k$.
(Our construction depends on the choice of the fixed point $u\in\d$.)

\myskip
The paper is organized as follows:
We start
in Section~\ref{prelim} with some general remarks on the Lie groups $\SU$ and $\tSU$
and their embeddings in the $4$-dimensional pseudo-Euclidean space resp.\ in a certain $\r_+$-bundle, the universal cover of a positive cone in that pseudo-Riemannian space.
We describe in section~\ref{elements} some elements of the construction, such as affine half-spaces and their substitutes in the $\r_+$-bundle.
We also define prismatic sets $\Qx$, certain finite intersections of half-spaces, and study their properties.
After that we are prepared to state in section~\ref{results} our main results, Theorems A and B, and to prove them.
In section~\ref{examples} we describe our explicit computations of fundamental domains for particular pairs of discrete subgroups and give pictures of these fundamental domains.

\myskip
I would like to thank Egbert Brieskorn and Ludwig Balke for useful conversations related to this work.

\section{Preliminaries}

\label{prelim}

\noindent
We consider the $4$-dimensional pseudo-Euclidean space $E^{2,2}$ of signature $(2,2)$.
We think of $E^{2,2}$ as the real vector space $\c^2\cong\r^4$ with the symmetric
bilinear form
$$\<(z_1,w_1),(z_2,w_2)\>=\Re(z_1\bar z_2-w_1\bar w_2).$$

\myskip\noindent
In the pseudo-Euclidean space $E^{2,2}$ we consider the quadric
\begin{align*}
 \Grp
  &=\left\{a\in E^{2,2}\st\<a,a\>=-1\right\}\\
  &=\left\{(z,w)\in E^{2,2}\st|z|^2-|w|^2=-1\right\}.
\end{align*}
For a fixed $z\in\c$ the intersection
$$\{w\in\c\st(z,w)\in\Grp\}=\{w\in\c\st|w|^2=|z|^2+1\}$$
is the circle of radius $\sqrt{|z|^2+1}\ge1$.
It holds $|w|\ge1$ for any $(z,w)\in\Grp$.
The bilinear form on $E^{2,2}$ induces a Lorentz metric of signature $(2,1)$ on $\Grp$.
The quadric~$\Grp$ is a model of the pseudo-hyperbolic space.

\myskip
Furthermore we consider the cone over $\Grp$
$$
  \Lrp
  =\r_+\cdot\Grp
  =\{\la\cdot a\st\la>0,~a\in\Grp\}.
$$
The cone $\Lrp$ can be described as
\begin{align*}
 \Lrp
  &=\left\{a\in E^{2,2}\st\<a,a\><0\right\}\\
  &=\left\{(z,w)\in E^{2,2}\st|z|<|w|\right\}.
\end{align*}
For a fixed $z\in\c$ the intersection
$$\{w\in\c\st(z,w)\in\Lrp\}=\{w\in\c\st|w|>|z|\}$$
is the complement of the disc of radius $|z|$.
It holds $w\ne0$ for any $(z,w)\in\Lrp$.
The bilinear form on $E^{2,2}$ induces a pseudo-Riemannian metric of signature~$(2,2)$ on $\Lrp$.

\myskip
We may think of $\Lrp$ as a $\r_+$-bundle over $\Grp$ with radial projection
$\bdlmap:\Lrp\to\Grp$ as bundle map.
The map $\Lrp\to\d$ defined by $(z,w)\mapsto z/w$ is a principal $\c^*$-bundle,
where the action of $\la\in\c^*$ is defined by $\la\cdot(z,w)=(\la^{-1}z,\la^{-1 }w)$.
Let $\covmap:\tGrp\to\Grp$ be the universal covering.
Henceforth we identify the Lie group $\SU$ with $\Grp$ via
$$\WZZW\mapsto(z,\bar w),$$
and $\tSU$ with $\tGrp$.
The biinvariant metrics on $\Grp$ and $\tGrp$
are proportional to the Killing forms.
We denote the pull-back $\tLrp\to\tGrp$
of the $\r_+$-bundle $\bdlmap:\Lrp\to\Grp$ under the covering map
$\covmap:\tGrp\to\Grp$ also by $\covbdlmap$.
The following diagram commutes
$$
  \begin{CD}
   \tLrp            @>\covmap>> \Lrp          \\
   @V{\covbdlmap}VV         @VV{\bdlmap}V \\
   \tGrp            @>\covmap>> \Grp          \\
  \end{CD}
$$
$\Grp$ resp.\ $\tGrp$ is canonically embedded in $\Lrp$ resp.\ $\tLrp$ and
therefore there exist canonical trivializations $\Lrp\cong\Grp\times\r_+$  resp.\ $\tLrp\cong\tGrp\times\r_+$.
The covering $\tLrp$ inherits canonically a pseudo-Riemannian metric from $\Lrp$.

\myskip
We now give a brief description of the full isometry group of $\tGrp$
(compare sections 2.1--2.3 in \cite{KR}).
The product $\tGrp\times\tGrp$ acts on $\tGrp$ via
$$(g,h)\cdot x=gxh^{-1}$$
by Lorentz isometries since the metric is biinvariant.
The identity component $\Isom_0(\tGrp)$ of the isometry group is isomorphic to $(\tGrp\times\tGrp)/\De_Z$,
where
$$\De_Z=\{(z,z)\st z\in Z\}$$
and $Z$ is the centre of $\tGrp$.
The full isometry group of $\tGrp$ has four components
corresponding to time- and/or space-reversals.
Let $\inv$ be the geodesic symmetry at the identity given by $g\mapsto g^{-1}$
and $\spieg$ the lift of the conjugation by the matrix
$\left(\begin{smallmatrix} 0&1 \\ 1&0 \end{smallmatrix}\right)$
in $\Grp$ fixing the identity.
Then $\inv$ preserves the space-orientation and reverses the time-orientation,
while $\spieg$ reverses both the space- and time-orientation.
Moreover, the group $\Isom^+(\tGrp)=\<\Isom_0(\tGrp),\spieg\>$
is the full group of orientation-preserving isometries and
$$
  \Isom(\tGrp)
  =\<\Isom_0(\tGrp),\spieg,\inv\>
  \cong\Isom_0(\tGrp)\semiprod(\<\spieg\>\times\<\inv\>)
$$
is the full isometry group of $\tGrp$.


\myskip
The universal covering $\covmap:\tLrp\to\Lrp$ of
\begin{align*}
  &\Lrp=\left\{(z,w)\in E^{2,2}\st|z|<|w|\right\} \\
  \intertext{can also be described as}
  &\tLrp=\{(z,\al,r)\in\c\times\r\times\r_+\st|z|<r\}, \\
  &\covmap(z,\al,r)=(z,re^{i\al}).
\end{align*}
We call the number $\al\in\r$ the argument of the element $(z,\al,r)\in\tLrp$.

\myskip
The restriction of the covering map $\covmap:\tLrp\to\Lrp$ gives the description
of the universal covering $\covmap:\tGrp\to\Grp$ of
\begin{align*}
  &\Grp=\left\{(z,w)\in E^{2,2}\st|z|^2-|w|^2=-1\right\} \\
  \intertext{as}
  &\tGrp=\{(z,\al,r)\in\c\times\r\times\r_+\st |z|^2=r^2-1\}, \\
  &\covmap(z,\al,r)=(z,re^{i\al}).
\end{align*}
For $(z,\al,r)\in\tGrp$ the positive real number $r$ can be computed from $z$
and $\al$, hence we can also identify $\tGrp$ with $\c\times\r$ via
$(z,\al,r)\mapsto(z,\al)$.

\myskip
The map $\covbdlmap:\tLrp\to\tGrp$ can be described as
$$
  \covbdlmap(z,\al,r)
  =\left(\la^{-1}z,\al,\la^{-1}r\right)
  \quad\hbox{with}\quad\la=\sqrt{r^2-|z|^2}.
$$

\section{The Elements of the Construction}

\label{elements}

\noindent
For $g\in\tGrp$ let $\Eg$ resp.\ $\Ig$ be the connected component
of $\covmap^{-1}(\bEbg)$ resp.\ $\covmap^{-1}(\bIbg)$ containing~$g$, where
$\barg:=\covmap(g)$ is the image of $g$ in $\Grp$,
\begin{align*}
  \bEbg&:=\{a\in\Lrp\st\<g,a\>=-1\}=T_{\barg}\Grp\cap\Lrp \\
  \intertext{is the intersection of the affine tangent space~$T_{\barg}\Grp$ on $\Grp$ in the point $\barg$ with~$\Lrp$ and}
  \bIbg&:=\{a\in\Lrp\st\<g,a\>\le-1\}=T_{\barg}^-\Grp\cap\Lrp
\end{align*}
is the intersection the half-space~$T_{\barg}^-\Grp$ of $\c^2$
bounded by $\bEbg$ and not containing~$0$ with~$\Lrp$.
$\bEbg$ and $\bIbg$ are simply connected and even contractible,
hence their pre-images under the covering map $\covmap$ consist of infinitely
many connected components, one of them containing~$g$.

\myskip
The three-dimensional submanifold~$\Eg$ subdivides~$\tLrp$ in two
connected components, the closure of one of them is~$\Ig$, and we denote the
closure of the other by~$\Hg$.
The boundary of~$\Ig$, resp.~$\Hg$, is equal to~$\Eg$.

\myskip
As an example, for the unit elements $e=(0,0,1)$ in~$\tGrp$ and $\bare=\covmap(e)=(0,1)$ in~$\Grp$,
we have
$$\bIbe=\{(z,w)\in\c^2\st\Re(w)\ge1,~|z|<|w|\},$$
the boundary $\bEbe$ of $\bIbe$ is a one-sheeted hyperboloid of revolution.
The pre-image of $\bIbe$ is
$$
  \covmap^{-1}(\bIbe)
  =\{(z,\al,r)\in\c\times\r\times\r_+\st r\cdot\cos\al\ge1,~|z|<r\}.
$$
The connected components of $\covmap^{-1}(\bIbe)$ resp.\ $\covmap^{-1}(\bEbe)$
containing~$e$ are
\begin{align*}
  \Ie
  &=\left\{
     (z,\al,r)\in\c\times\r\times\r_+
     \st
     |\al|<\frac{\pi}{2},~r\ge\frac{1}{\cos\al},~|z|<r
   \right\}\\
  \intertext{and}
  \Ee
  &=\left\{
     (z,\al,r)\in\c\times\r\times\r_+
     \st
     |\al|<\frac{\pi}{2},~r=\frac{1}{\cos\al},~|z|<r
   \right\}.
\end{align*}
The subsets $\Eg$ resp.\ $\Ig$ have the analogous properties
because $\Eg=g\cdot\Ee$ and $\Ig=g\cdot\Ie$.

\myskip
We make use of the following construction (compare \cite{Mi}).
Given a base-point $x\in\d$ and a real number $t$,
let $\rho_x(t)\in\PSU$ denote the rotation through angle~$t$ about the point~$x$.
Thus we obtain a homomorphism $\rho_x:\r\to\PSU$,
which clearly lifts to the unique homomorphism $r_x:\r\to\tSU$
into the universal covering group.
Since $\rho_x(2\pi)=\Id_{\d}$,
it follows that the lifted element $r_x(2\pi)$ belongs
to the central subgroup $Z$ of $\tSU$.
Note that this element $r_x(2\pi)\in Z$ depends continuously on $x$, and
therefore is independent of the choice of $x$.
We easily compute $r_0(2t)=(0,-t,1)$
and hence $r_x(2\pi)=r_0(2\pi)=(0,-\pi,1)$ for all~$x\in\d$.
Moreover we obtain
\begin{align*}
  &r_0(2t)\cdot(z,\al,r)=(ze^{it},\al-t,r),\\
  &(z,\al,r)\cdot r_0(2t)=(ze^{-it},\al-t,r),\\
  &(z,\al,r)\cdot r_0(-2t)=(ze^{it},\al+t,r).
\end{align*}

\myskip
Let $\G_1$ and~$\G_2$ be discrete subgroups of finite level~$k$ in~$\tSU$.
For~$i=1,2$, let $\bG_i$ be the image of~$\G_i$ in~$\PSU$.
We assume the existence of a joint fixed point $u\in\d$ of~$\bG_1$ and~$\bG_2$.

\myskip
For~$i=1,2$, the isotropy group~$(\bG_i)_u$ of $u$ in $\bG_i$ is a finite cyclic group generated by $\rho_u(2\pi/p_i)$, where $p_i=|(\bG_i)_u|$.
The isotropy group $(\G_i)_u$ of $u$ in $\G_i$ is an infinite cyclic group generated by $\deck_i:=r_u(2\thet_i)$, where~$\thet_i=\frac{\pi k}{p_i}$.
We can assume without loss of generality that $u=0\in\d$.
Under this assumption it follows
$$
  \deck_i=r_0(2\thet_i)=\left(0,-\thet_i,1\right)
  \quad\hbox{and}\quad
  \deck_i\cdot(z,\al,r)=\left(ze^{i\thet_i},\al-\thet_i,r\right).
$$

\myskip\noindent
Now let us start with the construction of fundamental domains for the action of~$\G_1\times\G_2$ on~$\tGrp$.
For a point $x$ in the orbit $\G_1(u)$ let $\Tx$ be
$$\Tx=\{(g_1,g_2)\in\G_1\times\G_2\st g_1(u)=x\}.$$
Let
$$\Qx=\bigcap\limits_{(g_1,g_2)\in\Tx}H_{g_1g_2}.$$

\myskip\noindent
As an example, for $x=u$ we have that
$$\Tu=(\G_1)_u\times\G_2=\{(\deck_1^{m_1},\deck_2^{m_2})\st m_1,m_2\in\z\}=\<(d_1,e),(e,\deck_2)\>.$$
The generator~$(\deck_1,e)$ acts on $\tGrp$ by left multiplication
$$\deck_1\cdot(z,\al,r) = (ze^{i\thet_1},\al-\thet_1,r).$$
The generator~$(e,\deck_2)$ acts on $\tGrp$ by right multiplication
$$(z,\al,r)\cdot\deck_2^{-1} = (ze^{i\thet_2},\al+\thet_2,r).$$

\myskip\noindent
Let $p=\lcm(p_1,p_2)$ be the least common multiple of~$p_1$ and~$p_2$.
Let
$$\deck=r_u(2\pi k/p)=r_u(2\thet),\quad\text{where}~\thet=\frac{\pi k}{p}.$$
The element~$\deck$ acts on $\tGrp$ by left multiplication
$$\deck\cdot(z,\al,r) = (ze^{i\thet},\al-\thet,r)$$
and it acts on the $(\al,r)$-half-plane by the translation mapping
$$\tau(\al,r)=(\al-\thet,r).$$
An important assumption for the following construction is
$$p>k.$$
In terms of the element~$\deck$ the assumption~$p>k$ means that the argument $\thet$ of $\deck$ is less then $\pi$.

\myskip\noindent
We have
$$
  \Qu
  =\bigcap\limits_{(g_1,g_2)\in T(u)}H_{g_1g_2}
  =\bigcap\limits_{m_1,m_2\in\z}H_{\deck_1^{m_1}\deck_2^{m_2}}
  =\bigcap\limits_{m\in\z}H_{\deck^m},
$$
since~$\<d_1,d_2\>=\<d\>$.

\begin{figure}
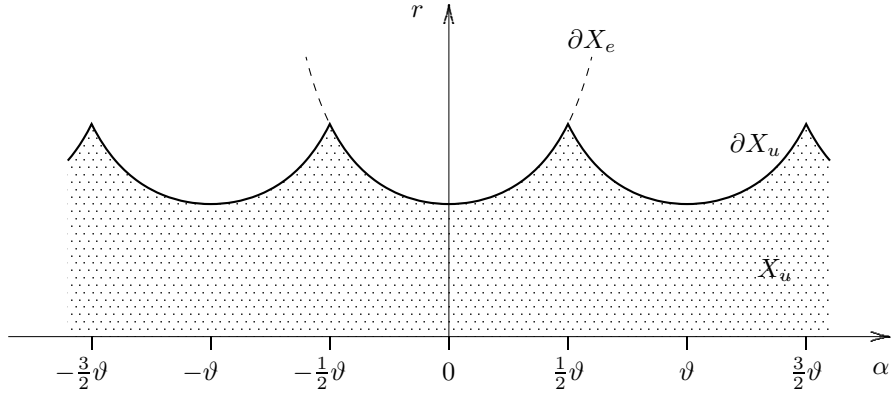

  \begin{center}
    \forcehmode
    \Pictexure{.}{qu}
  \end{center}
  \caption {The image $X_u$ of $Q_u$ in the $(\alpha,r)$-half-plane}
  \label {figa}
\end{figure}

\myskip\noindent
What does the set
$$\Qu=\bigcap\limits_{m\in\z}H_{d^m}$$
look like?
The image of the set~$\He$ under the projection $(z,\al,r)\mapsto(\al,r)$ is
$$\Xe=\{(\al,r)\in\r\times\r_+\st r\cdot\cos\al\le1\enskip\text{or}\enskip|\alpha|\ge{\pi}/{2}\}.$$
The images of the sets $H_{\deck^m}=\deck^m\cdot\He$ under the projection $(z,\al,r)\mapsto(\al,r)$
are the translates $\tau^m(\Xe)$ of the set~$\Xe$.
The manifold $\Qu$ is a disc bundle over its image
$\Xu=\bigcap_{m\in\z}\tau^m(\Xe)$ in the $(\al,r)$-plane.
The shaded area in figure~\ref{figa} is $\Xu$.
(The real line is not part of~$\Xu$.)
The subsets~$\Qx$ are images of the subset~$\Qu$ under the action of the group~$\G_1\times\G_2$.
For any~$x\in\G_1(u)$ there is an element~$g\in\G_1$ such that~$g(x)=u$.
Then~$\Qx=g\cdot\Qu$.


The manifolds $g\Qu$ play a central role in our construction.
We want to explain the geometric nature of these objects.
We have described $\Qu$ as a disc bundle over the set $\Xu$ in the $(\al,r)$-half-plane $\r\times\r_+$.
We may describe $\Qu\subset\tLrp\subset\c\times\r\times\r_+$ as
$$\Qu=(\c\times\Xu)\cap\tLrp.$$
We think of $\Xu$ as a universal covering of a punctured plane polygon.
Consider the following diagram of covering maps
$$
  \begin{tabular}{*2{c@{}}c}
    $\r\times\r_+$
          & \cdarrowheightfalse\cdarrowxx{\pi'}{}{0}{1}
          & $\c^\text{\rlap{$*$}}$ \\
    & \cdarrowxx{}{\pi~}{-45}{1.4142} & \!\cdarrowv{}{\pi''}{-90}{1} \\
    &                                   & $\c^\text{\rlap{$*$}}$
  \end{tabular}
$$
where $\pi(\al,r)=re^{i\al}$, $\pi'(\al,r)=r^{1/k}e^{i\al/k}$ and $\pi''(z) = z^k$.
We now consider the curve~$\pi(\dd\Xu)$.
It is easy to see that this is a regular star polygon $\big\{\frac{2p}{k}\big\}$ when $k$ is odd and a regular star polygon $\big\{\frac{p}{k}\big\}$ when $k$ is even,
whereby a star polygon~$\big\{\frac{n}{m}\big\}$, with~$n$ and~$m$ positive integers, is a figure formed by connecting with straight lines every $m$-th point
out of $n$~regularly spaced points lying on a circle (see H.S.M.~Coxeter~\cite{Co}, \S2.8, pp. 36--38).

\myskip\noindent
{\bf Remark:} $k=2$, $p_1=5$, $p_2=3$, $p=15$: $\big\{\frac{15}{2}\big\}$~star polygon.

\myskip\noindent
Therefore the curve $\pi'(\dd\Xu)$ is a curvilinear $2p$-gon
covering the star polygon once or twice.
Let $P'\subset\c$ and $P=\Pu\subset\c$ be the plane areas bounded by the curvilinear polygon
$\pi'(\dd\Xu)$ and by the star polygon $\pi(\Xu)$. The images of
$\Xu$ are the punctured plane polygons $\pi'(\Xu)=P'\setminus\{0\}$ and
$\pi(\Xu)=P\setminus\{0\}$.
We think of the product $\c\times P'$ as a 4-dimensional $2p$-gonal {\it prism}.
$\c\times\Xu$ is the universal covering of the pierced prism $\c\times(P'\setminus\{0\})$.
The product $\c\times P\subset\c^2$ might be considered as a 4-dimensional {\it ``star prism''.}
Its axis $\c\times\{0\}$ does not meet $\Lrp\subset\c\times\c^*$.
Therefore the universal covering $\pi:\tLrp\rightarrow\Lrp$ maps $\Qu$
to the intersection of $\Lrp$ with the star prism:
$$\pi(\Qu)=\Lrp\cap(\c\times\Pu).$$

\myskip
In the following lemma we prove some properties of the sets $\Qx$.
We first give some definitions.
Let $s:\tGrp\to\r_+$ be a section
in the bundle~$\tLrp\cong\tGrp\times\r_+$.
We call the set
$$\{(a,\la)\in\tGrp\times\r_+\st \la=s(a)\}$$
the {\it graph}\/ of~$s$ and the set
$$\{(a,\la)\in\tGrp\times\r_+\st \la\le s(a)\}$$
the {\it subgraph}\/ of~$s$.

\begin{lem}
\label{Qx}
For a point $x\in\d$ in the orbit~$\G_1(u)$ of the point~$u$
under the action of the group~$\G_1$ the following holds:
\begin{enumerate}[(i)]
\item
For any point $(z,w)\in\covmap(\Qx)$
$$|w|-|z|\le|w-\bar xz|\le f(|x|),$$
where
$$f(t):=\frac{\sqrt{1-t^2}}{\cos\frac{\thet}{2}}.$$
\item
The set $\Qx$ is a subgraph of a section
in the bundle $\tLrp\cong\tGrp\times\r_+$,
while its boundary is the graph of this section.
\end{enumerate}
\end{lem}

\begin{proof}
Our proof is in two steps.
We first check the properties of $\Qx$ in the case $x=u$.
In this case the properties follow from the explicit description of the set $\Qu$.
Then we use the fact that for any $x\in\G_1(u)$ there is an element $g\in\G_1$ such that $\Qu=g\cdot\Qx$ to prove the properties of $\Qx$ for $x\ne u$.

\myskip
Let us first describe explicitly the image $\Xu$ of the set $\Qu$ in the $(\al,r)$-plane $\r\times\r_+$.
The set $\Xu$ is the shaded area in figure~\ref{figa}.
It is a subgraph of a function $\r\to\r_+$.
Let us denote this function by~$\varphi$.
We now describe the function~$\varphi$ explicitly.
The function~$\varphi$ is periodic with period~$\thet$,
hence it is sufficient to describe~$\varphi$ on $[-\thet/2,\thet/2]$.
For $\al\in[-\thet/2,\thet/2]$ it holds
$$\varphi(\al)=\frac{1}{\cos\al}.$$
For any $\al\in\r$ it holds
$$\varphi(\al)\le\frac{1}{\cos\frac{\thet}{2}}$$
(with equality for $\al=(2k+1)\thet/2$, $k\in\z$).

\myskip
Now let us verify the first assertion of the lemma.
The inequality
$$|w|-|z|\le|w-\bar xz|$$
follows from $|z|<|w|$ and $|x|<1$.
It remains to prove the second inequality.

\myskip
Let us verify the first assertion of the lemma in the case $x=u$.
(Recall that we assumed $u=0$.)
For $x=u=0$ the second inequality in the first part of the lemma reduces to
$$|w|\le\frac{1}{\cos\frac{\thet}{2}}$$
for any point $(z,w)\in\covmap(\Qu)$.
Let us consider a point $(z,w)\in\covmap(\Qu)$ and its preimage $(z,\al,r)\in\Qu$.
By definition of the map $\covmap$ it holds $w=r e^{i\al}$.
For the point $(z,\al,r)\in\Qu$ it holds $(\al,r)\in\Xu$.
The set~$\Xu$ is the subgraph of the function~$\varphi$, hence
$$r\le\varphi(\al)\le\frac{1}{\cos\frac{\thet}{2}}$$
for any point $(\al,r)\in\Xu$.
Hence
$$|w|=r\le\frac{1}{\cos\frac{\thet}{2}}.$$

\myskip
Let us verify the first assertion of the lemma for any $x$.
Let us consider a point $x\in\G_1(u)$ and an element $g\in\G_1$ such that $g(x)=u$.
Let $(a,b)\in\Grp$ be the image of the element~$g$ under~$\covmap$.
The element $(a,b)\in\Grp$ corresponds to the matrix
$$\begin{pmatrix} \bar b&a \\ \bar a&b \end{pmatrix}\in\SU$$
and acts on~$\d$ by
$$(a,b)\cdot x=\frac{\bar b x+a}{\bar a x+b}.$$
The property $(a,b)\cdot x=u=0$ implies $a=-\bar bx$.
From $(a,b)\in\Grp$ we conclude
$$-1=|a|^2-|b|^2=|-\bar b x|^2-|b|^2=-|b|^2\cdot(1-|x|^2)$$
and hence
$$|b|={1\over\sqrt{1-|x|^2}}.$$
Let us consider $(z,w)\in\covmap(\Qx)$ and $(z',w')=g\cdot(z,w)\in\covmap(\Qu)$.
On the one hand $(z',w')\in\covmap(\Qu)$ implies
$$|w'|\le\frac{1}{\cos\frac{\thet}{2}}.$$
On the other hand
$$
  |w'|
  =|\bar az+b w|
  =|-b\bar xz+bw|
  ={1\over\sqrt{1-|x|^2}}\cdot|w-\bar x z|.
$$
Hence
$$|w-\bar xz|\le\frac{\sqrt{1-|x|^2}}{\cos\frac{\thet}{2}}.$$

\myskip
Let us verify the second assertion of the lemma in the case $x=u$.
For the set~$\Qu$ we can describe the corresponding section $s_u:\tGrp\to\r_+$ explicitly as
$$s_u(z,\al,r)=\frac{\varphi(\al)}{r}.$$

\myskip
Let us verify the second assertion of the lemma for any $x$.
Let us consider a point $x\in\G(u)$ and an element $g\in\G$ such that $\Qu=g\cdot\Qx$.
Then the section $s_x:\tGrp\to\r_+$ is given by
$$s_x(a)=s_u(g\cdot a).\qedhere$$
\end{proof}

\begin{lem}
\label{Qxlocfin}
The family $(\Qx)_{x\in\G_1(u)}$ is locally finite
in the sense that any point of $\tLrp$ has a neighbourhood intersecting
only finitely many prisms $\Qx$.
\end{lem}


\begin{proof}
We prove that the family $(\covmap(\Qx))_{x\in\G_1(u)}$ is locally finite (in $\Lrp$).
This fact implies the local finiteness of the family $(\Qx)_{x\in\G(u)}$,
since if a subset $U$ of $\Lrp$ has an empty intersection with $\covmap(\Qx)$ then the intersection of the pre-image $\covmap^{-1}(U)$ with $\Qx$ is empty too.
By lemma~\ref{Qx}(i) for any point $x\in\G_1(u)$ and any point $(z,w)\in\covmap(\Qx)$ the difference $|w|-|z|$ is bounded from above by~$f(|x|)$.
The values~$f(t)$ tend to zero as~$t$ tends to~$1$.
Choosing a point $(z_0,w_0)\in\Lrp$ and a positive number $\ve<|w_0|-|z_0|$,
the neighbourhood $U:=\{(w,z)\in\Lrp\st|w|-|z|>\ve\}$ of the point $(z_0,w_0)$ can intersect $\covmap(\Qx)$ only for $|x|$ sufficiently small (so that $f(|x|)>\ve$).
But the group $\G_1$ is discrete, so there are only finitely many points~$x$ in~$\G(u)$ with norm~$|x|$ under a given bound.
This finishes the proof.
\end{proof}

\begin{rem}
This property of $\Qx$ allows us to deal
with $P=\cup\Qx$ in a similar way as with a finite union of polytopes.
\end{rem}

\begin{lem}
\label{EgQxlocfin}
The family $(\Eg\cap Q_{g(u)})_{g\in\G_1}$ is locally finite.
\end{lem}

\begin{proof}
This is immediate from the local finiteness of the family $(\Qx)_{x\in\G_1(u)}$
plus the easy observation that the family $(\Eg\cap Q_{g(u)})_{g\in(\G_1)_u}$ is locally finite.
\end{proof}

We consider in $\tLrp$ the four-dimensional polytope
$$
  P:=\bigcup\limits_{x\in\G_1(u)}\Qx
  =\bigcup\limits_{x\in\G_1(u)}\bigcap\limits_{g\in\Tx}\Hg.
$$

\begin{lem}
\label{ddPtGrp-homeo}
The projection~$\ddP\to\tGrp$ is a $\G_1\times\G_2$-equivariant homeomorphism.
\end{lem}


\begin{proof}
From lemma~\ref{Qx}(ii) we know that the set $\Qx$ is a subgraph of a section in the bundle $\tLrp\cong\tGrp\times\r_+$.
A union of a locally finite family of subgraphs of sections in~$\tLrp$ is again a subgraph of a section in~$\tLrp$.
To see this, let us first consider the following toy version of this statement:
A union of subgraphs of functions~$f_1,\dots,f_k:\r\to\r_+$ is again a subgraph of a function~$f:\r\to\r_+$, where~$f=\max(f_1,\dots,f_k)$.
This is clear in the toy case and generalizes to the case of a locally finite family of subgraphs of sections in~$\tLrp$.
Thus the polyhedron $P=\cup\Qx$ is a subgraph of a section in the bundle~$\tLrp\cong\tGrp\times\r_+$ as a union of a locally finite family of subgraphs.
But for a subgraph of a section in the bundle $\tLrp$
it is clear that the bundle map $\tLrp\to\tGrp$ induces a homeomorphism from its boundary (equal to the graph of the section) onto $\tGrp$.
This homeomorphism is $\G_1\times\G_2$-equivariant since the projection $\tLrp\to\tGrp$ is $\G_1\times\G_2$-equivariant.
\end{proof}

\section{The Main Results}

\label{results}

\noindent
Now we can state the main result

\begin{thmA}
The boundary of~$P$ is invariant with respect to the action of~$\G_1\times\G_2$.
The subset
$$\Fg=\Cl_{\ddP}(\Int(\dd\Hg\cap\ddP))$$
is a fundamental domain for the action of~$\G_1\times\G_2$ on~$\ddP$.
The family
$$(F_{g_1g_2})_{g_1\in\G_1, g_2\in\G_2}$$
is locally finite in~$\ddP$.
The projection~$\tLrp\to\tGrp$ induces a $\G_1\times\G_2$-equivariant homeomorphism
$$\ddP\to\tGrp.$$
The image~$\calFg$ of~$\Fg$ under the projection is a fundamental domain for the action of~$\G_1\times\G_2$ on~$\tGrp$.
The family~$(\calF_{g_1g_2})_{g_1\in\G_1, g_2\in\G_2}$ is locally finite.
For every elements $g_1,h_1\in\G_1$, $g_2,h_2\in\G_2$ with~$g_1g_2\ne h_1h_2$
the intersection~$\calF_{g_1g_2}\cap\calF_{h_1h_2}$ lies in a totally geodesic submanifold of~$\tGrp$.
\end{thmA}

\begin{rem}
In this section all closures are taken in~$\ddP$.
We use the shorthand $\Cl$ instead of $\Cl_{\ddP}$.
\end{rem}

\begin{lem}
\label{lem-gen-top}
Let $X$ be a topological space.
Let $A$ and $B$ be closed subsets of $X$.
Then
\begin{enumerate}[{\rm (i)}]
\item
\label{gen-top-intclintcl=intcl}
$\Int\Cl\Int A=\Int A$,
\item
\label{gen-top-intcl-and-clintcl}
$\Int A\cap\Cl\Int B\ne\emptyset\Rightarrow\Int(A\cap B)\ne\emptyset$.
\end{enumerate}
\end{lem}

\begin{lem}
\label{lem-about-Fg}
$$\Int\Fg=\Int(\Eg\cap\ddP)\quad\text{and}\quad\Cl\Int\Fg=\Fg.$$
\end{lem}

\begin{proof}
The assertions follow from Lemma~\ref{lem-gen-top}(\ref{gen-top-intclintcl=intcl}) with $A=\Eg\cap\dd P$.
\end{proof}

\begin{proof}
To prove that $\Fg$ is a fundamental domain we have to prove two properties.
The first property is that the images of $\Fg$ have no common inner points,
\ie the intersection $\Int(\Fg)\cap\Fh$ is empty if $g\ne h$.
The second property is that $\Cl(\cup_{g\in\G}\Int\Fg)=\ddP$,
\ie roughly speaking the images of $\Fg$ cover the whole space~$\ddP$.

\myskip
Let us first prove that the intersection $\Int(\Fg)\cap\Fh$ is empty if $g\ne h$.
Suppose on the contrary that there are elements~$g,h\in\G$
such that $g\ne h$ and $\Int(\Fg)\cap\Fh\ne\emptyset$.
Let us consider the closed subsets $A=\Eg\cap\dd P$ and $B=\Eh\cap\dd P$.
By Lemma~\ref{lem-about-Fg} it holds $\Int(\Fg)=\Int A$,
hence the assumption $\Int(\Fg)\cap\Fh\ne\emptyset$ can be rewritten as $\Int A\cap\Cl\Int B\ne\emptyset$.
From Lemma~\ref{lem-gen-top}(\ref{gen-top-intcl-and-clintcl}) it follows that $\Int(A\cap B)\ne\emptyset$.
This means that the set $\Int(\Eg\cap\Eh\cap\ddP)$ is not empty.
But since the totally geodesic submanifolds~$\Eg$ and~$\Eh$ intersect transversally,
the intersection $\Eg\cap\Eh$ has no inner points in~$\ddP$.

\myskip
Since $\Fg\subset\Eg\cap Q_{g(u)}$ lemma~\ref{EgQxlocfin} implies that
the family~$(\Fg)_{g\in\G}$ is locally finite in~$\ddP$.
Lemma~\ref{ddPtGrp-homeo} says that the projection~$\ddP\to\tGrp$ is a $\G$-equivariant homeomorphism.

\myskip
Now let us prove the property $\Cl(\cup_{g\in\G}\Int\Fg)=\ddP$.
Since
$$
  \Cl\big(\bigcup\limits_{g\in\G}\Int\Fg\big)
  \supset\bigcup\limits_{g\in\G}\Cl\Int\Fg
  =\bigcup\limits_{g\in\G}\Fg
$$
(where the last equality holds by Lemma~\ref{lem-about-Fg}),
it suffices to prove that $\cup_{g\in\G}\Fg=\ddP$.
Consider $a\in\ddP$.
From the definition of~$P$ and local finiteness (according to Lemma~\ref{EgQxlocfin})
of the family $(\Eg\cap Q_{g(u)})_{g\in\G}$ it follows that
in some neighbourhood of the point~$a$ only finitely many elements of~$\G$ are relevant,
\ie there exists a neighbourhood~$U$ of the point~$a$ in~$\tLrp$ and elements $g_1,\dots,g_n\in\G$ such that
$$\ddP\cap U=\bigcup\limits_{i=1}^n\,(E_{g_i}\cap\ddP\cap U).$$
We may assume without loss of generality that
the map $\covmap|_U:U\to\covmap(U)$ is a homeomorphism.
The image of $P\cap U$ under this homeomorphism is an intersection of an open
subset of~$\Lrp$ with a finite union of finite intersections of
half-spaces~$\Hg$ with the property $a\in\dd\Hg$.
Suppose that $a\not\in\Cl\Int(E_{g_i}\cap\ddP)=F_{g_i}$ for all $i\in\{1,\dots,n\}$.
This is only possible if for each~$i\in\{1,\dots,n\}$
the set $E_{g_i}\cap\ddP\cap U$ is contained in a $2$-dimensional submanifold of~$\tLrp$.
Thus $\ddP\cap U$ is contained in the union of finitely many $2$-dimensional submanifolds.
On the other hand it follows from lemma~\ref{ddPtGrp-homeo}
that $\ddP$ is homeomorphic to a $3$-dimensional manifold~$\tGrp$.
This contradiction implies that $a\in\Fg$ for some~$g\in\G$.
\end{proof}

\begin{lem}
\label{lem-desr-ddP}
The boundary $\ddP$ of $P=\cup_{x\in\G(u)}\Qx$ can be described as follows
$$
  \ddP
  =\dd\big(\bigcup\limits_{x\in\G(u)}\Qx\big)
  =\big(\bigcup\limits_{x\in\G(u)}\dd\Qx\big)\bs\big(\bigcup\limits_{x\in\G(u)}\Int\Qx\big).
$$
This means that a point $p$ is in the boundary of $P$
if and only if $p$ is not an interior point of any $\Qx$ with $x\in\G(u)$
and $p$ is a boundary point of $\Qx$ for some $x\in\G(u)$.
\end{lem}

\begin{proof}
From lemma~\ref{Qx}(ii) we know that the set $\Qx$ is a subgraph of a section~$s_x$ in the bundle $\tLrp\cong\tGrp\times\r_+$
$$\Qx=\{(a,\la)\in\tGrp\times\r_+\st\la\le s_x(a)\}.$$
The set~$P=\cup\Qx$ is the subgraph of the section $s_P=\max s_x$.
(In this proof $\max$ means $\max_{x\in\G(u)}$, $\cup$ means $\cup_{x\in\G(u)}$, $\exists\,x$ means $\exists\,x\in\G(u)$ and so on.)
This property would be obvious for a finite union of subgraphs.
Using local finiteness (according to Lemma~\ref{Qxlocfin}) we prove that this property also holds for~$P$.
But for a subgraph
$$X=\{(a,\la)\in\tGrp\times\r_+\st\la\le s(a)\}$$
of a section~$s$ in the bundle $\tLrp$ it is clear that $(a,\la)\in\dd X$ if and only if $\la=s(a)$.
Hence
$$(a,\la)\in\ddP\iff\la=s_P(a).$$
By definition of $s_P$
$$\la=s_P(a)\iff(\exists\,x\hskip10pt\la=s_x(a))\quad\text{and}\quad(\forall\,x\hskip10pt \la\ge s_x(a)).$$
On the other hand
\begin{align*}
  (a,\la)\in\cup\,\dd\Qx&\iff\exists\,x\hskip10pt \la=s_x(a),\\
  (a,\la)\not\in\cup\Int\Qx&\iff\forall\,x\hskip10pt \la\ge s_x(a).
\end{align*}
\end{proof}

\begin{lem}
\label{IntFe-in-ddQu}
$\Int\Fe\subset\dd\Qu$.
\end{lem}

\begin{proof}
By Lemma~\ref{lem-about-Fg} it holds $\Int\Fe=\Int(\Ee\cap\ddP)$.
Suppose that there is a point $a\in\Int\Fe=\Int(\Ee\cap\ddP)$ such that $a\not\in\dd\Qu$.
Since $a\in\dd P$ and $a\not\in\dd\Qu$ there exists $x\in\Guou$ such that $a\in\dd\Qx$.
Then any neighbourhood of $a$ intersects $\Ee\cap\Int\Qx\subset\Ee\bs\ddP$.
The projection $\covbdlmap:\tLrp\to\tGrp$ is continuous
and the restriction $\covbdlmap|_{\ddP}:\ddP\to\tGrp$ is a homeomorphism,
therefore any neighbourhood of $a$ intersects
$((\covbdlmap|_{\ddP})^{-1}\circ\covbdlmap)(\Ee\bs\ddP))\subset\ddP\bs\Ee$.
This implies $a\not\in\Int(\Ee\cap\ddP)=\Int\Fe$. Contradiction.
\end{proof}

\begin{prop}
\label{FG=Fe}
$$\Fe=\Cl\Int\left((\Ee\cap\dd\Qu)-(\cupIntQxGuou)\right).$$
\end{prop}

\begin{proof}
Let $\hF:=(\Ee\cap\dd\Qu)-(\cup_{x\in\Guou}\Int\Qx)$.
We claim that $\Fe$ and $\hF$ coincide up to the boundary, \ie $\Int\Fe=\Int\hF$.
To prove this we show the inclusions $\Int\Fe\subset\Int\hF$ and $\Int\hF\subset\Int\Fe$.
We first prove that $\Int\Fe\subset\Int\hF$.
To that end we show that $\Int\Fe\subset\hF$.
Then $\Int\Int\Fe\subset\Int\hF$ and $\Int\Fe=\Int\Int\Fe$ imply $\Int\Fe\subset\Int\hF$.
To see that $\Int\Fe$ is contained in $\hF$ we have to show (by definition of~$\hF$)
that $\Int\Fe$ is contained in~$\Ee$, in~$\dd\Qu$, and does not intersect~$\Int\Qx$
for all~$x\in\Guou$.
By definition of~$\Fe$ it holds $\Int\Fe\subset\Ee$.
By Lemma~\ref{IntFe-in-ddQu} it holds $\Int\Fe\subset\dd\Qu$.
Finally for any $x\in\Guou$ it holds $\Fe\cap\Int\Qx=\emptyset$
because of the fact that $\Fe$ is contained in $\ddP$, and
$\ddP\cap\Int\Qx=\emptyset$ by Lemma~\ref{lem-desr-ddP}.
This implies $\Int\Fe\subset\hF$ and therefore $\Int\Fe\subset\Int\hF$.
We now have to prove the inclusion $\Int\hF\subset\Int\Fe$.
From the definition of~$\hF$ it follows that $\hF\subset\Ee$.
Moreover $\hF\subset\dd\Qu\subset(\cup_{x\in\Gu}\dd\Qx)$ and
$\hF\cap(\cup_{x\in\Guou}\Int\Qx)=\emptyset$
imply by Lemma~\ref{lem-desr-ddP} that $\hF\subset\ddP$.
Now from $\hF\subset\Ee\cap\ddP$ it follows
that $\Int\hF\subset\Int(\Ee\cap\ddP)=\Int\Fe$,
where the last equality holds by Lemma~\ref{lem-about-Fg}.
We now have proved both inclusions, \ie we know that $\Int\hF=\Int\Fe$.
From this it follows that $\Cl\Int\hF=\Cl\Int\Fe=\Fe$.
\end{proof}

\begin{lem}
\label{Fe-compact}
If $\G$ is co-compact, then $\Fg$ is compact.
\end{lem}

\begin{proof}
Consider a sequence $a_k$ in~$\Int\Fg$.
Let $\varphi$ be the composition of the projection maps $\ddP\to\tGrp$ and $\tGrp\to\tGrp/\G$.
Since the quotient~$\tGrp/\G$ is compact we may assume without loss of
generality that the sequence $\varphi(a_k)$ tends to a limit $\bar a\in\tGrp/\G$.
Since~$\varphi$ is surjective there exists a pre-image~$a\in\ddP$ of~$\bar a$ under~$\varphi$.
Hence there is a sequence $h_k$ in~$\G$ such that the sequence $h_k a_k$ tends to~$a$.
Since the family $(\Fg)_{g\in\G}$ is locally finite there exists a
neighbourhood~$U$ of~$a$ that intersects only finitely many fundamental domains~$\Fg$.
Therefore the set $\{h_k|k\in\n\}$ is finite.
After choosing a subsequence we may assume that the sequence $h_k$ is constant, say $h_k=h$.
Then the sequence $h a_k$ tends to~$a$,
hence the sequence $a_k$ tends to $h^{-1}a$.
This implies $h^{-1}a\in\Fg$.
\end{proof}

\begin{thmB}
If $\G$ is co-compact then  $\Fg$ is a compact polyhedron,
i.e. a finite union of finite compact intersections of half-spaces~$\Ia$.
\end{thmB}

\begin{proof}
The family $(\Qx)_{x\in\G(u)}$ is locally finite
and the fundamental domain $\Fe$ is compact by lemma~\ref{Fe-compact}.
From this it follows that there is a finite subset $E\subset\Gu$ such that
$\Fe\cap\Qx=\emptyset$ for all $x\in\Gu\bs E$.
By proposition~\ref{IntFe-in-ddQu} this implies the assertion.
\end{proof}

\section{Examples}

\label{examples}

\myskip
We have computed the fundamental domains explicitly for those infinite series of pairs of discrete subgroups which correspond via the construction described in~\cite{Pr:qgor}
to certain series of $\q$-Gorenstein quasi-homogeneous surface singularities.
In particular the quotient of $\tSU$ by one of the corresponding group action is diffeomorphic to the link of the corresponding quasi-homogeneous singularity.

\myskip
A discrete co-compact subgroup $\Gamma$ of level $k$ in $\tSU$ such that the image in $\PSU$ is a triangle group with signature $(\alpha_1,\alpha_2,\alpha_3)$
will be denoted by~$\Gamma(\alpha_1,\alpha_2,\alpha_3)^k$.

\myskip
The following figures show some of the explicitly computed fundamental domains.

\myskip
Some explanations are required to make the figures of fundamental domains comprehensible.
The image $\pi(F_e)$ of the fundamental domain $F_e$ is a compact polyhedron in~$\su(1,1)$ with flat faces.
The Lie algebra~$\su(1,1)$ is a $3$-dimensional flat Lorentz space of signature $(n_+,n_-) = (2,1)$.
Such a polyhedron has a distinguished rotational axis of symmetry.
The direction of this axis is negative definite, and the orthogonal complement is positive definite.
Changing the sign of the pseudo-metric in the direction of the rotational axis transforms Lorentz space into a well-defined Euclidean space.
The image~$\pi(F_e)$ of the fundamental domain is then transformed into a polyhedron in Euclidean space with dihedral symmetry.
Figures~\ref{fund-533-3-2}, \ref{fund-733-3-2} and \ref{fund-933-3-2} show the Euclidean polyhedra obtained in this way
in the cases~$\G(5,3,3)^2\times(C_3)^2$, $\G(7,3,3)^2\times(C_3)^2$ and~$\G(9,3,3)^2\times(C_3)^2$.
The direction of the rotational axis is vertical.
The top and bottom faces are removed.

\myskip
The polyhedra in figures~\ref{fund-533-3-2}, \ref{fund-733-3-2} and~\ref{fund-933-3-2} are all scaled by the same factor
to illustrate the proportions between different fundamental domains.

\myskip
Figures~\ref{ident533-3-2}, \ref{ident733-3-2}, \ref{ident933-3-2-a}, \ref{ident933-3-2-b} illustrate the identification schemes
for the cases~$\G(5,3,3)^2\times(C_3)^2$, $\G(7,3,3)^2\times(C_3)^2$ and~$\G(9,3,3)^2\times(C_3)^2$.
The face identification is equivariant with respect to the dihedral symmetry of the polyhedron.
The faces labeled with the same letter and shaded in the same way are identified.
Numbers on the edges of shaded faces indicate the identified flags (face, edge, vertex).

\vskip1in

\input fund-533-3-2.pic

\input ident-533-3-2.pic

\input fund-733-3-2.pic

\input ident-733-3-2.pic

\input fund-933-3-2.pic

\input ident-933-3-2-a.pic

\input ident-933-3-2-b.pic